\documentclass{amsart}
\usepackage{amssymb, amsthm, amsmath, amsfonts, graphicx}
\input xy
\xyoption{all}

\begin{document}

\newtheorem {lemma} {LEMMA} [section]
\newtheorem {thm}[lemma]{THEOREM}
\newtheorem {prop}[lemma]{PROPOSITION}
\newtheorem {cor}[lemma]{COROLLARY}
\newtheorem {defn}[lemma]{DEFINITION}

\numberwithin{equation}{section}

\newcommand{\pf}[1]{\textbf{Proof #1:}}
\newcommand{\CZ}{(\mathbb C^*)^{\mathbb Z}}
\newcommand{\C}{\mathbb C}
\newcommand{\ignore}[1]{}
\newcommand {\BOX} {\rule{2mm}{2mm}}

\def\ds{\displaystyle}

    \title{The solenoid and holomorphic motions for H\'enon maps}             
    \author{Philip P. Mummert}                                
    \address{ Department of Mathematics\\
      Purdue University\\
      West Lafayette, IN 47907-2067}                                    
    \email{pmummert@math.purdue.edu}                                      
    \urladdr{}
    \date{04/18/2007}                                       
    \begin{abstract} The McMullen-Sullivan holomorphic motion for topologically conjugate, complex polynomials with connected Julia set follows level sets of the B\"ottcher coordinate.  The Buzzard-Verma holomorphic motion for hyperbolic, unstably connected, polynomial diffeomorphisms of $\C^2$ follows level sets of the Bedford-Smillie solenoid map.  It follows that this solenoid map is a conjugacy for those H\'enon maps that are perturbations of (one-dimensional) hyperbolic maps with connected Julia set.  \end{abstract} 
\maketitle

\section{Introduction}

In this paper we show the equivalence of two holomorphic motions. One motion is derived from the external rays and equipotential lines for the connected Julia set of a complex polynomial.  The other is obtained by canonically extending the natural motion on the grand orbits of the critical points of the polynomial.  An analogous scenario plays out in two complex variables.

  From well-known results of dynamics in one complex variable, for a polynomial map $P$ on $\C$ of degree $d\ge 2$, there are two possibilities: Either the filled Julia set $K$ (equivalently, the set of all points in $\C$ that are bounded under iteration of $P$) contains all the critical points of $P$, in which case $K$ (and $J=\partial K$) is connected, or alternatively, if at least one critical point lies outside $K$, then $K$ is disconnected.  When $J$ is connected (which, in the case of quadratic polynomials, exactly corresponds to those parameter values that lie in the famed Mandelbrot set), there is an important conjugacy known as the B\"ottcher coordinate, a holomorphic bijection $\phi: \C \setminus K \to \C \setminus \mathbb D$, where $\mathbb D$ denotes the closed unit disk.  The dynamics of $P$ are conjugated by $\phi$ to the action of $z \mapsto z^d$ on the exterior of the closed unit disk: $\phi P(z)=(\phi(z))^d$. A simple dynamical model indeed.  Note that the B\"ottcher coordinate equips $\C \setminus K$ with polar coordinates via $\phi^{-1}$. Lines of constant angle are mapped via $\phi^{-1}$ to ``external rays'' on $\C \setminus K$. Much study has been made of the combinatorics of the landing and pinching of these external rays. For more details on one-dimensional complex dynamics see \cite{M}.

  McMullen and Sullivan \cite{MS} prove the following: For a holomorphic family of topologically conjugate complex polynomials $P_\lambda$, $\lambda \in \Delta$, there is a holomorphic motion $\psi_\lambda$ on $\C$ that respects the dynamics, i.e. $P_\lambda \psi_\lambda=\psi_\lambda P_0$. 

  We point out that this motion preserves the B\"ottcher coordinate.  Such a result is part of the folklore on deformation theory of polynomial maps, i.e. the B\"ottcher coordinate is a holomorphic conjugacy invariant under deformation. However, we present a complete statement of proof here.  Since $\psi_\lambda$ is quasiconformal on $\C \setminus K_0$, the result follows from the canonical nature of the construction of the motion rather than the uniqueness of the B\"ottcher coordinate. The following diagram commutes:

\[\xymatrix{
& \C \setminus K_\lambda \ar[rr]^{P_\lambda} \ar[dd]|\hole_(.3){\phi_\lambda} & & \C \setminus K_\lambda \ar[dd]_{\phi_\lambda}\\
\C \setminus K_0 \ar[rr]^(.65){P_0} \ar[ur]^{\psi_\lambda} \ar[dr]_{\phi_0} & & \C \setminus K_0 \ar[ur]^{\psi_\lambda} \ar[dr]_{\phi_\lambda} & \\
& \C \setminus \mathbb D \ar[rr]^{z\mapsto z^d} && \C \setminus \mathbb D \\
}\]

\begin{thm} Given a family of topologically conjugate polynomials with connected Julia set, $P_\lambda$ for $\lambda \in \Delta$ with the McMullen-Sullivan holomorphic motion $\psi_\lambda(z)$ defined on $\hat \C$ over $\Delta$ and the B\"ottcher coordinate $\phi_\lambda$, then \[ \phi_0(z)=\phi_\lambda \psi_\lambda(z) \] for all $\lambda \in \Delta$ and $z\in \C \setminus K_0$. In particular, $\psi_\lambda$ is holomorphic on $\C \setminus K_0$ for each fixed $\lambda \in \Delta$.
\end{thm}

Moving to higher dimensions, any dynamically interesting polynomial diffeomorphism $H$ of $\C^2$ is in fact conjugate to a composition of H\'enon mappings, i.e. maps of the form \[(x,y)\mapsto (y, P(y)-bx)\] where $b \in \C \setminus \{0\}$, and $P(y)$ is complex polynomial of degree $\ge 2$ (See \cite{BuS} for an expository introduction). Following \cite{BV}, for brevity we use the generalized term ``H\'enon map'' to describe a nonelementary polynomial diffeomorphism of $\C^2$.  Let $K^+$ be the set of points bounded under iteration of $H$ and $K^-$ those bounded under iteration of $H^{-1}$. Let $J^+=\partial K^+$, $J^-=\partial K^-$, and the Julia set $J=J^+ \cap J^-$.  Hubbard and Oberste-Vorth point out the existence of an analytic function $\varphi^+$, analogous to $\phi$, defined for $|y|$ large and $|y|\ge |x|$ in $\C^2 \setminus K^+$, again mapping to $\C \setminus \mathbb D$ such that $\varphi^+ H (x,y) = (\varphi^+(x,y))^d$, where $d= {\rm deg}\; H$.  However, $\varphi^+$ cannot be analytically continued to all of $\C^2 \setminus K^+$. See \cite{HO1}.
 
Drawing similarities to the dynamics of one complex variable, there is also an analogous notion for the connectivity of $J$; when $J$ is ``unstably connected'' Bedford and Smillie show that $\varphi^+$ analytically extends to a neighborhood of $J^-_+ = J^- \setminus K^+$.  See \cite{BS6}.

Because $H$ is invertible, any conjugate model must be an invertible dynamical system.  The degree-$d$ exterior solenoid $\Sigma_+ \subset \C^\mathbb Z$ is given by the inverse limit of $z \mapsto z^d$ on $\C\setminus \mathbb D$. Let $p=(x,y)\in \C^2$.  When $J$ is unstably connected,  the map $\Phi: J^-_+ \to \Sigma_+$, where $\Phi(p)_t = \varphi^+(H^t (p))$ for each $t\in \mathbb Z$, semiconjugates the dynamics of $H$ to the shift map $\sigma$ on the exterior solenoid. That is, $\Phi H(p) = \sigma \Phi(p)$.  $\Phi$ is known to be a finite-degree covering map  \cite{BS7}.  If $\Phi$ is invertible, then this is, in fact, a conjugacy and the exterior solenoid proves to be a simple topological model of the dynamics of hyperbolic, unstably connected H\'enon maps. We show that to be the case for those hyperbolic, unstably connected H\'enon maps that are ``perturbations'' of one-dimensional maps.

Like the B\"ottcher coordinate in one complex variable, $\Phi$ equips $J^-_+$ with external rays.  If $\Phi$ is injective, then these rays in fact provide a unique coding for points in $J$.  Note that Bedford and Smillie circumvent questions about the injectivity of $\Phi$ in  \cite{BS7} section 4, by carefully crafting a conjugating homeomorphism that roughly corresponds to taking a $k$-th root of $\Phi$, where $k$ is the covering degree of $\Phi$.  This procedure allows the formulation of their results on solenoidal external rays.

Analogous to $\psi_\lambda$,  Buzzard and Verma define a holomorphic motion that respects the dynamics for hyperbolic H\'enon maps  \cite{BV}.  In particular, if $\Delta$ is a complex polydisk corresponding to a region of hyperbolic parameter values, there is a conjugating holomorphic motion $\Psi_\lambda$ for $\lambda \in \Delta$ mapping the action of $H_0$ on $J^-_+(0)$ to that of $H_\lambda$ on $J^-_+(\lambda)$.  That is, $\Psi_\lambda H_0 = H_\lambda \Psi_\lambda$.  Here $J^-_+(\lambda)$ denotes the set $J^-_+$ for the map $H_\lambda$.

Our first question is whether this holomorphic motion is compatible with the function $\Phi$ mentioned above, similar to the one variable case summarized in Theorem 1.1.  The following diagram commutes:

\[\xymatrix{
& J^-_+(\lambda) \ar[rr]^{H_\lambda} \ar[dd]|\hole_(.3){\Phi_\lambda} & & J^-_+(\lambda) \ar[dd]_{\Phi_\lambda}\\
J^-_+(0) \ar[rr]^(.65){H_0} \ar[ur]^{\Psi_\lambda} \ar[dr]_{\Phi_0} & & J^-_+(0) \ar[ur]^{\Psi_\lambda} \ar[dr]_{\Phi_\lambda} & \\
& \Sigma_+ \ar[rr]^{\sigma} && \Sigma_+\\
}\]

\begin{thm} Given a family of hyperbolic, unstably connected H\'enon maps, $H_\lambda$ for $\lambda \in \Delta$ with the Buzzard-Verma holomorphic motion $\Psi_\lambda$ defined on $U^+_0$ over $\Delta$ and the Bedford-Smillie solenoidal map $\Phi_\lambda$, then \[ \Phi_0(p) = \Phi_\lambda \Psi_\lambda(p)\] for all $\lambda \in \Delta$ and $p \in J^-_+(0)$. In particular, $\Psi_\lambda$ is holomorphic on each leaf of $J^-_+(0)$ for each fixed $\lambda \in \Delta$.
\end{thm}

\begin{cor} The covering degree of $\Phi_\lambda$ is constant for $\lambda \in \Delta$.  
Also, the identifications of external rays are preserved for  $\lambda \in \Delta$.
\end{cor}

In fact, a similar result holds for H\'enon maps with Jacobian parameter $b \to 0$.  If $P$ is a hyperbolic polynomial with connected Julia set, consider the set of $b\in \C\setminus \{0\}$ such that the H\'enon map, $H_{b}(x,y)=(y, P(y)-bx)$ is hyperbolic and unstably connected.  Let $\mathcal M_P$ be a component of this set such that the closure of $\mathcal M_P$ contains $0$.

By a well-known result of Hubbard and Oberste-Vorth  \cite{HO2}, $\mathcal M_P$ contains the set $\{b: 0<|b|<B\}$ for some $B>0$ depending on $P$.

\begin{thm} If $b\in \mathcal M_P$, then the Bedford-Smillie solenoidal map $\Phi_b$ is a topological conjugacy from $J^-_+(b)$ to $\Sigma_+$. \end{thm}

The proofs of Theorems 1.1 and 1.2 appear in section 6; the proof of Theorem 1.4 appears in section 7.
I am tremendously grateful to Greg Buzzard for our countless discussions and his helpful explanations.  I also thank John Smillie for many fruitful conversations. Lastly, thanks to Kevin Pilgrim for pointing out some corrections. 

\section{Preliminaries}

Consider the family of diffeomorphisms of $\C ^2$ known as the complex H\'enon mappings, with $P_j$ a complex polynomial of degree $\ge 2$, $x,y,b_j \in \C$ and $b_j\neq 0$:
\[F_j(x,y)= (y, P_j(y)-b_jx)\]

We generalize the term ``H\'enon map'' to include any non-elementary, polynomial diffeomorphism $H$ of $\C^2$, since these are conjugate to a composition $F_1 \circ \hdots \circ F_m$ where $F_j$ is as above (see \cite{FM}).
Let $K^+$ denote those points $p=(x,y)\in \C^2$ bounded under forward iteration of $H$. $J^+$ denotes $\partial K^+$.
Similarly define $K^-$ and $J^-$ for backward iteration of $H$, that is, iteration under $H^{-1}$.
Let the Julia set $J= J^- \cap J^+$.  Write $U^+= \C^2 \setminus K^+$ and $U^-= \C^2 \setminus K^-$. The set $J^-_+=J^- \cap U^+$ will be of particular interest. 

We refer to $b=\prod_{j=1}^m b_j$ as the Jacobian parameter, since $b=\det DH$.  Replacing $H$ by $H^{-1}$ if necessary, assume that $|b|\le 1$.  From  \cite{BS1} and  \cite{BS2}, if $H$ is hyperbolic, (that is, complementary directions of uniform expansion and contraction on $J$), the unstable set of $J$, $W^u(J)$, is $J^-$ minus finitely many periodic sinks.  The stable set of $J$, $W^s(J) = J^+$.  The interior of $K^+$ is the sink basins, and $K^-$ has empty interior.  $W^u(J)$ has a dynamically defined lamination whose leaves are the unstable manifolds of points in $J$. Each unstable manifold leaf is conformally equivalent to $\C$.  Each leaf may intersect $U^+$ in a number of components; this forms a lamination of $J^-_+$.

Bedford and Smillie \cite{BS6} say a H\'enon map is ``unstably connected with respect to a saddle point $p$'' if some component of $W^u(p) \cap U^+$ is simply connected.  This is in fact equivalent to each component of $W^u(q) \cap U^+$ being simply connected for any saddle periodic point $q$.  In this case we say that $H$ is ``unstably connected.''  Furthermore, as Bedford and Smillie prove, if $|b|\le 1$, the condition that $H$ is unstably connected is equivalent to $J$ being connected.

Let $\lambda$ be an element of a closed (for compactness arguments that follow) complex polydisk, $\overline \Delta$, that corresponds to a set of parameter values that give rise to a holomorphic family of degree-$d$ H\'enon maps $H_\lambda$ that are hyperbolic and unstably connected.  We will often refer to the parameter value $0\in \Delta$ as a basepoint.  We will also abuse notation throughout, writing $\Delta$ to refer always to a region of parameter space; sometimes the open disk in $\C$ and other times the open polydisk in $\C^n$, depending upon the context.

Define the following regions in $\C^2$:

$\ds V^+(R) = \{(x,y)\in \C^2: |y| \ge |x|, |y| \ge R \}$

$\ds V^-(R) = \{(x,y)\in \C^2: |x| \ge |y|, |x| \ge R \}$

$\ds V(R) = \{(x,y)\in \C^2: |x| < R, |y| < R \}$

From  \cite{BS1} Lemma 2.1, there exists an $R_\lambda$ such that these partitioning regions $V^+, V^-$, and $V$ satisfy various properties for $R \ge R_\lambda$:

\begin{enumerate}
\item $V^+ \subset U^+$, $\ds U^+ = \bigcup_{n\ge 0} H^{-n}(V^+)$
\item $V^- \subset U^-$, $\ds U^- = \bigcup_{n\ge 0} H^{n}(V^-)$
\item $H(V^+) \subset V^+$ and $H^{-1}(V^-) \subset V^-$
\item $H(V)\subset V \cup V^+$ and $H^{-1}(V) \subset V \cup V^-$
\end{enumerate}
Also, $K^+ \subset V \cup V^+$, $K^- \subset V \cup V^-$.
Tracking the forward iterates of any point: it either enters and remains in $V^+$, or, if the point lies in $K^+$, it enters and remains in $V$.  

We will often consider the collection of all these sets over a region of parameter space, rather than focusing on one individual parameter value (as has been the usual setting in the papers of Hubbard and Oberste-Vorth and Bedford and Smillie).  Subscripts will be written as needed to indicate reference to a particular parameter value.  Let $\ds R_\Delta$ satisfy the necessary conditions above for every $\lambda \in \Delta$. Define the following sets in $\Delta \times \C^2$:

$\mathcal V = \Delta \times V(R_\Delta)$

$\mathcal V^+ = \Delta \times V^{+}(R_\Delta)$

$ \mathcal U^{+} = \{(\lambda, p): \lambda \in \Delta, p \in U^{+}_\lambda\}$

$ \mathcal J^{+} = \{(\lambda, p): \lambda \in \Delta, p \in J^{+}_\lambda\}$

Similarly, define the sets $\mathcal V^-$, $\mathcal U^-$, and $\mathcal J^-$ and let $\mathcal J^-_+ := \mathcal J^- \cap \mathcal U^+$.

Observe $\mathcal U^{+}$ is an open set. ($U^{+}_{\lambda_0}$ is open for any $\lambda_0$ and given $p \in U^{+}_{\lambda_0}$, $ p \in U^{+}_\lambda$ for $\lambda$ sufficiently close to $\lambda_0$).

In order to keep the notion of iteration and a dynamical system on these sets, let $\mathcal H: \Delta \times \C^2 \to \Delta \times \C^2$ be given by \[\mathcal H(\lambda, p)=(\lambda, H_\lambda(p))\] 

The collection of laminations of $J^-_+(\lambda)$ for each $\lambda \in \Delta$ forms a ``leafwise trivial holomorphic family of laminations'' on $\mathcal J^-_+$.  These terms are explicitly defined in  \cite{BV}.

In $\mathcal V^+$, there is a choice of a $d^n$th root of $\pi_y \circ \mathcal H^n$ that approximates $y$.  In this way define 
$\varphi^+: \mathcal V^+ \to \C\setminus \mathbb D$ by \[\ds \varphi^+_\lambda (p):= \varphi^+ (\lambda, p) = \lim_{n\to \infty} \left ( \pi_y \circ H_\lambda^n(p) \right )^{1/d^n}\]

\begin{prop}[Hubbard-Oberste-Vorth] $\varphi^+$ is holomorphic on $\mathcal V^+$. \end{prop}

\pf{} Hubbard and Oberste-Vorth show $\varphi^+_\lambda$ is holomorphic in $p$ for each $\lambda$. Their proof rewrites $\varphi^+_\lambda$ as the limit of a telescoping infinite product.  For every $\lambda$, we have a uniform estimate on the product (independent of $p$ and $\lambda$), hence the convergence is locally uniform.  As the locally uniform limit of functions holomorphic in $\lambda$ (as well as $p$), the limit, $\varphi^+$, is holomorphic. See \cite{HO1} for more details on $\varphi^+$. \BOX \\

Let $\ds G^+: \Delta \times \C^2 \to \mathbb R$ be given by
\[  G^+_\lambda(p):= G^+(\lambda, p) = \lim_{n \to \infty} \frac{1}{d^n} \log_+ \|H_\lambda^n(p) \| \]

On $\mathcal V^+$, $G^+ = \log |\varphi^+|$.
The function $G^+$ is pluriharmonic on $\mathcal U^+$ by \cite{BS1} Prop. 3.3.  Hence for each $\lambda$, the complex 1-form $\partial G^+_\lambda$ is holomorphic on $U^+_\lambda$ and determines a holomorphic foliation $\mathcal G^+_\lambda$ on $U^+_\lambda$.
$G_\lambda^+(p)$ is pluriharmonic on $\mathcal U^+$, meaning it is locally the real part of a holomorphic function on $\mathcal U^+$.  The level sets of this holomorphic function correspond to the plaques of $\mathcal G_\lambda^+$, hence $\mathcal G_\lambda^+$, $\lambda \in \Delta$ is a holomorphic family of laminations.  In fact, $\mathcal G_\lambda^+$ is a leafwise trivial holomorphic family of laminations and biholomorphic to $\Delta \times \C$.  See \cite{BV}.

\section{Holomorphic motions}

For an introduction to holomorphic motions in one variable see \cite{AM}.  For holomorphic motions of several complex variables see \cite{J}.
\begin{defn} A {\bf holomorphic motion} of a set $E\subset \hat \C$ over $\Delta$ centered at $\delta_0 \in \Delta$ is a mapping $\alpha: \Delta \times E \to \hat \C$, such that \begin{enumerate}
\item $\alpha(\delta_0, z)=z$ for every $z \in E$.
\item $\alpha(\lambda_0, z)$ is injective on $E$ for each fixed $\lambda_0 \in \Delta$.
\item $\alpha(\lambda, z_0)$ is holomorphic in $\lambda$ for each fixed $z_0\in E$.
\end{enumerate}
\end{defn}

We will assume our motions are centered at $\delta_0=0$, unless otherwise indicated. For reference, we recall the theorem of Ma\~n\'e, Sad, and Sullivan \cite{MSS}, also known as the $\lambda$-Lemma, as well as a theorem of Bers and Royden \cite{BR}, also known as the Harmonic $\lambda$-Lemma.

\begin{thm}[Ma\~n\'e-Sad-Sullivan]\label{mss}  If $\alpha:\Delta \times E \to \hat \C$ is a holomorphic motion, then there is a unique extension to a holomorphic motion $\alpha: \Delta \times \overline E \to \hat \C$.  Furthermore, $\alpha$ is continuous in $(\lambda, z)$ and $\alpha(\lambda_0, z)$ is quasiconformal on $\overline E$ for each fixed $\lambda_0 \in \Delta$.
\end{thm}

\begin{defn} A Beltrami coefficient $\mu$ is {\bf harmonic} on an open set $U \subset \hat \C$ if \[\mu(z)=\frac{{\overline {\nu(z)}}}{\rho^2(z)}\] for $z\in U$ where $\rho^2(z)|dz|$ is the area element of the hyperbolic metric on $U$, and the function $\nu(z)$ is holomorphic in $z\in U$. \end{defn}

\begin{thm}[Bers-Royden]\label{br} If $\alpha:\Delta \times E \to \hat \C$ is a holomorphic motion, then there is a canonical extension to a holomorphic motion $\alpha: \Delta_{1/3} \times \hat \C \to \hat \C$ uniquely characterized by the following property:  $\alpha_\lambda$ has harmonic Beltrami coefficient on $\hat \C \setminus \overline E$ for each $\lambda \in  \Delta_{1/3}$.
\end{thm}

Next we introduce the holomorphic motions of Theorem 1.1 and 1.2 respectively, both constructed using Theorem 3.4:

\begin{thm}[McMullen-Sullivan 7.4] Let $P_0$ be a postcritically stable complex polynomial, or equivalently, there is a family of topologically conjugate polynomials $P_\lambda$ for $\lambda \in \Delta$.  Then there exists $r>0$ and a map \[\psi: \Delta_r \times \hat \C  \to \hat \C\] that is a holomorphic motion that respects the dynamics.  In particular, setting $\psi_\lambda(z):=\psi(\lambda, z)$, we have
\begin{enumerate}
\item $\psi_0(z)=z$.
\item $\psi_\lambda$ is a homeomorphism for each fixed $\lambda$.
\item $\psi_\lambda(z)$ is holomorphic in $\lambda$ for each fixed $z\in \hat \C$.
\item $\psi_\lambda P_0 = P_\lambda \psi_\lambda$.
\item $\psi_\lambda$ is quasiconformal with harmonic Beltrami coefficient on $\C \setminus K_0$ for each fixed $\lambda$.
\end{enumerate}
\end{thm}

The last property follows from the construction of the holomorphic motion. McMullen and Sullivan point out the existence of a unique motion respecting the dynamics on the grand orbits of the critical points of $P_\lambda$. Then the Bers-Royden theorem gives extension to all of $\hat \C$.  In particular, the extended motion has harmonic Beltrami coefficient on $\C \setminus K_0$.

\begin{thm}[Buzzard-Verma 1.2]  Let $H_\lambda$ be a one-parameter family of hyperbolic, unstably connected H\'enon maps depending holomorphically on $\lambda \in \Delta$. Then there exists $r>0$ and a map \[\Psi:\Delta_r \times (J^+_0 \cup U^+_0 \cup J^-_0) \to J^+_\lambda \cup U^+_\lambda \cup J^-_\lambda \]
such that defining $\Psi_\lambda(p):=\Psi(\lambda, p)$, we have
\begin{enumerate}
\item $\Psi_0(p)=p$.
\item $\Psi_\lambda$ is a homeomorphism for each fixed $\lambda$.
\item $\Psi_\lambda(p)$ is holomorphic in $\lambda$ for each fixed $p \in J^+_0 \cup U^+_0 \cup J^-_0$.
\item $\Psi_\lambda$ maps each leaf of $J_0^-$ to a leaf of $J_\lambda^-$. (Similarly with leaves of $J^+_0 \cup U^+_0$).
\item $\Psi_\lambda H_0 = H_\lambda \Psi_\lambda$ on  $J^+_0 \cup U^+_0 \cup J^-_0$.
\item $\tau_\lambda^{-1}\Psi_\lambda\tau_0: \C \to \C$ is a holomorphic motion with harmonic Beltrami coefficient on $\tau_0^{-1}(W^u(p_0) \setminus J_0)$,
\end{enumerate}  
where $p_0\in J_0$ and $\tau_\lambda:\C \to W^u(p_0)$ is the parametrization obtained using affine structures in \cite{BV} Theorem 5.6 such that $\tau$ is holomorphic in $(\lambda, z)$, $\tau_\lambda(0)=\Psi_\lambda(p_0)$, and $\tau_\lambda(1)=\Psi_\lambda(q_0)$ for some $q_0\in W^u(p_0)\setminus \{p_0\}$.
\end{thm}

Again, the use of Theorem 3.4 by Buzzard and Verma to construct their holomorphic motion gives the last property.  In particular, since $H_\lambda$ is unstably connected, $\tau_\lambda^{-1}\Psi_\lambda\tau_0$ has a harmonic Beltrami coefficient on the component $\tau_0^{-1} (W^u(p_0) \cap U^+_0)$ for every $\lambda\in \Delta$.

Let $\Psi^{-1}: \mathcal J^-_+ \to J^-_+(0)$ be given by $\Psi^{-1}(\lambda, p) = q$ where $\Psi(\lambda, q) = p$. Since $\Psi_\lambda$ is a homeomorphism, this map is well-defined and surjective.  Write $\Psi^{-1}_\lambda(p):=\Psi^{-1}(\lambda, p)$. 

\begin{prop} $\Psi^{-1}$ is continuous. \end{prop}

\pf{}
The proof of \cite{BV} Theorem 1.1 verifies that $\Psi_\lambda^{-1}$ is continuous on $J^-_+(\lambda)$ by showing $\Psi_\lambda$ is proper for each $\lambda \in \Delta$.  The proof here is identical: the map from $\Delta \times J^-_+(0)$ to $\mathcal J^-_+$ given by $(\lambda, p_0) \mapsto (\lambda, \Psi_\lambda(p_0))$ must also be proper, therefore $\Psi^{-1}$ is continuous. \BOX \\

\section{Extension of $\varphi^+$}

We extend some results about $\varphi^+_{\lambda_0}$ for fixed $\lambda_0$ to $\varphi^+$ viewed as a function of $\lambda$ as well as $p$, for all $\lambda$ in a region of parameter space $\Delta$.

\begin{defn} Let $A$ be a closed subset of a topological space $X$.  $A$ is a {\bf strong deformation retract} of $X$ if there is a continuous map $F: X \times [0,1] \to X$ such that for every $x\in X$, $a \in A$ and $t\in [0,1]$:
\begin{enumerate}
\item $F(x,0)=x$
\item $F(x,1)\in A$
\item $F(a,t)=a$
\end{enumerate}
\end{defn}
Thus $F$ is a homotopy between the identity map on $X$ and a retraction of $X$ onto $A$. 

\begin{lemma} [\cite{BS6} Lemma 2.5] Let $A \subset X$ be a strong deformation retract.  Let $\pi:\tilde{Y} \to Y$ be a covering map.  Let $\rho:X \to Y$ be a continuous map and assume we are given a continuous map $\sigma: A \to \tilde{Y}$ which is a lift of $\rho |_A$.  Then there is a unique continuous lift $\tilde{\sigma}$ of $\rho$ which agrees with $\sigma$ on $A$.
\end{lemma}

To say  $\tilde{\sigma}$ is a lift of $\rho$, means $\tilde{\sigma}$ maps $X$ to $\tilde{Y}$ continuously such that $\pi \circ \tilde{\sigma}=\rho$.  
The proof follows from the homotopy lifting property of the covering map $\pi$, the uniqueness of lifts of paths, and the uniqueness property of lifts of homotopies.

Bedford and Smillie use this lemma to prove:

\begin{lemma} [\cite{BS6} Theorem 2.6, Theorem 6.3] If $H_\lambda$ is unstably connected for a fixed $\lambda$ then $\varphi^+_\lambda$ has a unique continuous extension from $\ds J^-_\lambda \cap V^+_\lambda$ to $\ds J^-_\lambda \cap U^+_\lambda$ which satisfies the functional equation \[\varphi^+_\lambda(H_\lambda(p)) = (\varphi^+_\lambda(p))^d\]
In addition, for fixed $\lambda$, $\varphi^+_\lambda$ is holomorphic on leaves of $J^-_+(\lambda)$.
\end{lemma}

By applying an additional deformation step we obtain the analogous result:

\begin{thm} Given a family of hyperbolic and unstably connected H\'enon maps $H_\lambda$ depending holomorphically on $\lambda \in \Delta$, then $\varphi^+ = \varphi^+(\lambda, p)$ has a unique continuous extension from $\mathcal V^+$ to $\ds \mathcal V^+ \cup \mathcal J^-_+$ which satisfies the functional equation \[\varphi^+(\mathcal H (\lambda,p)) = (\varphi^+(\lambda,p))^d\]
\end{thm}

 For fixed $\lambda = 0$, Lemma 4.3 gives a continuous extension of $\varphi^+_0$ on $J^-_+(0)$.  Now there is a strong deformation retract from $\mathcal J^-_+$ to the slice at $J^-_+(0)$ by the Buzzard-Verma holomorphic motion.  So $\varphi^+_0$ uniquely extends to all of $\mathcal J^-_+$.  We must also check that this extension agrees with the one obtained using Lemma 4.3 for $\varphi^+_\lambda$ with fixed non-zero $\lambda$.  The details of the proof follow.

\pf{}
Let $V' \subset V^+(R_\Delta)$ be chosen such that $\Psi_\lambda(V') \subset V^+(R_\Delta)$ and \[\bigcup_{j \ge 0} H^{-j}_\lambda(V') \supset  J^-_\lambda \cap (U^+_\lambda \setminus V^+(R_\Delta))\] for every $\lambda \in \Delta$.
For instance, take $V'$ to be a set in $V^+(R_\Delta)$ which is relatively compact and contains a fundamental domain of $J^-_\lambda$ for each $\lambda \in \Delta$.  By the continuity of $\Psi$ and the compactness of $\overline \Delta$ such a neighborhood exists.

Let  $\mathcal V' := \{(\lambda,\Psi_\lambda(V'): \lambda \in \Delta\} \subset \mathcal V^+$. Fix $n\in \mathbb N_0$. We wish to apply Lemma 4.2 and will name our spaces and maps accordingly. Let  $\pi:\C \setminus \mathbb D \to \C \setminus \mathbb D$ be given by $\pi(z)=z^{d^n}$. Set $\tilde{Y}$ and $Y$ equal to $\C \setminus \mathbb D$. Then $\pi:\tilde{Y} \to Y$ is a degree-$d^n$ covering map. Let 
\[ A_{\lambda} := J^-_\lambda \cap \bigcup_{j=0}^n H^{-j}_\lambda\left (\Psi_\lambda(V')\right ) \]
\[ \mathcal A := \{(0, p_0): p_0 \in A_{0} \} \]
\[ X :=\{(\lambda,p): \lambda \in \Delta, p \in A_{\lambda}\} \]
Write $A_\lambda^0, \mathcal A^0$, and $X^0$ for the above sets with $n=0$.   

 $\mathcal A$ is a closed subset of $X$ (the slice at $\lambda=0$), $X \subset \mathcal J^-_+$, and $X^0 \subset \mathcal V^+$.
Also $X^0 \subset X$, $A_\lambda^0 \subset A_\lambda$ and $\mathcal A^0 \subset \mathcal A$ for any $n$. Let $\rho:X \to Y$ be given by $\rho := \varphi^+ \circ \mathcal H^n$.
$\rho$ is holomorphic since $ \mathcal H^n(X) \subset \mathcal V^+$. 
  
Lemma 4.3 gives a continuous extension of $\varphi^+_0$ from $A_{0}^0$ to $A_{0} \subset J^-_0 $ which defines the map $\sigma$  on $\mathcal A$ by $\sigma(0,p_0) = \varphi^+_0(p_0)$.
$\sigma$ is indeed a lift of $\ds \rho |_{\mathcal A}$, since for $(0,p_0) \in \mathcal A$,
$\rho (0,p_0)= \varphi^+_0 \circ H^n_0(p_0)=(\varphi^+(0,p_0))^{d^n}$.

Let
\[F(\lambda, p, t)=\Big (\lambda(1-t), \Psi \big (\lambda(1-t), \Psi^{-1}(\lambda, p) \big )\Big )\]

  $F$ is a strong deformation retract from $X$ to $\mathcal A$ as:
\begin{enumerate} 

\item
$F(\lambda, p, 0) = \big (\lambda, \Psi_\lambda(\Psi^{-1}_\lambda(p))\big ) = (\lambda,p)$ for $(\lambda,p) \in X$.
\item
$F(X, 1) = (0,\Psi^{-1}(X)) =  \Big (0, J^-_0 \cap  \bigcup_{j=0}^n H^{-j}_\lambda\left (\Psi_\lambda(V')\right ) \Big )  \subset \mathcal A$
\item
$F(0, p_0, t) = \big (0, \Psi_0(\Psi^{-1}_0(p_0))\big ) = (0,p_0)$ for $(0,p_0) \in \mathcal A$
\end{enumerate}

By Lemma 4.2, there is a unique continuous map $\tilde{\sigma}$ on $X$ that agrees with $\sigma$ on $\mathcal A$ and satisfies $\rho=\pi \circ \tilde{\sigma}$. So for $(0, p) \in \mathcal A \subset X$, $\tilde{\sigma}(0,p)= \sigma(0, p)= \varphi^+_0(p)$ and for  $(\lambda, p)\in X$, $\varphi^+ \mathcal H^n(\lambda, p)  = (\tilde{\sigma}(\lambda, p))^{d^n}$.

It remains to show $\tilde{\sigma}=\varphi^+$ on $X$, not just on $\mathcal A$.  This is not immediate from the uniqueness of the lift $\tilde{\sigma}$ since we do not know a priori that $\varphi^+$ is continuous in the variable $\lambda$ on all of $X$.  However, $\varphi^+$ is continuous on $X^0$ and by the uniqueness of Lemma 4.2, one can show that for any $n$, $\tilde{\sigma}=\varphi^+$ on $X^0$.


Next, fix $\lambda_0 \in \Delta$. 
Let  $\ds \tilde{\sigma}_{\lambda_0}: A_{\lambda_0}  \to \C \setminus \mathbb D$ be given by $\tilde{\sigma}_{\lambda_0}(p):= \tilde{\sigma}(\lambda_0, p)$.
We have shown that $ \tilde{\sigma}_{\lambda_0} = \varphi^+_{\lambda_0} $ on $A^0_{\lambda_0}$. 
Now, for $p \in A_{\lambda_0}$
\[\varphi^+(\mathcal H^n (\lambda_0,p)) = (\tilde{\sigma}(\lambda_0,p))^{d^n}\]
\[\tilde{\sigma}_{\lambda_0}(H^n_{\lambda_0}(p)) = (\tilde{\sigma}_{\lambda_0}(p))^{d^n}\]
since $H^n_{\lambda_0}(p) \in A^0_{\lambda_0}$.
Again by the uniqueness of the lift in Lemma 4.2, any continuation from $A^0_{\lambda_0}$ to  $A_{\lambda_0}$ satisfying the above functional equation must be unique. Therefore, $\tilde{\sigma}_{\lambda_0} = \varphi^+_{\lambda_0}$ on $A_{\lambda_0}$.  
As this can be done for every $\lambda_0 \in \Delta$, it follows that  $\tilde{\sigma} =  \varphi^+$ on $X$.
Since $\ds \mathcal J^-_+ \subset \mathcal V^+ \cup \bigcup_{n\ge 0} X_n$, conclude that $\varphi^+$ has a unique continuous extension to $\mathcal J^-_+$ which satisfies the functional equation $\varphi^+ \circ \mathcal H^n  = (\varphi^+)^{d^n}$. \BOX

\begin{lemma}[\cite{BS6} Lemma 6.2] If $H_\lambda$ is unstably connected, then $\varphi^+_\lambda$ has an analytic continuation to a neighborhood of $J^-_+$. \end{lemma}

The same method of proof as that of Bedford and Smillie \cite{BS6} (Lemma 6.2) yields the analogous result:

\begin{lemma} If $H_\lambda$ is a family of hyperbolic and unstably connected H\'enon maps depending holomorphically on $\lambda \in \Delta$, then $\varphi^+ = \varphi^+(\lambda, p)$ has an analytic continuation to a neighborhood of $\mathcal J^-_+$. \end{lemma}

\pf{}  
For any point in the holomorphic family of laminations $\mathcal G^+_\lambda$, there is a holomorphic family of leaves through it, $L_\lambda$.  Let $(\lambda_1, p_1)$ and $(\lambda_2, p_2)$ be in the same family of leaves, that is $p_1$ lies in $L_{\lambda_1}$, a leaf of the lamination $G^+(\lambda_1)$ and  $p_2$ lies in $L_{\lambda_2}$, the corresponding leaf of the lamination $G^+(\lambda_2)$.
There exists a biholomorphism between $\{(\lambda, p) : \lambda \in \Delta, p \in L_\lambda\}$ and $\Delta \times \C$, since this family of leaves is leafwise trivial (see \cite{BV}). For some choice of such a biholomorphism, this induces a metric on the family of leaves. Let $d_{\mathcal G}((\lambda_1, p_1),(\lambda_2, p_2))$ denote the distance with respect to the induced metric.

For $(\lambda_0, p_0) \in \mathcal U^+$, consider the ``nearest neighbors'': those points $(\lambda, p) \in \mathcal J^-_+$ which also lie in the family of leaves of $\mathcal G^+_\lambda$ through $(\lambda_0, p_0)$ and minimize the $d_{\mathcal G} $ distance among all such points. 

Let $P$ be an element of $\mathcal J^-_+$.  We want to show $\varphi^+$ has an analytic continuation to a neighborhood of $P$.  
By the leafwise triviality of the lamination and by the transversality of the lamination of $\mathcal J^-_+$ with the lamination $\mathcal G_\lambda^+$, one can choose local coordinates $\mu , u$, and $v$ so that $B:=\{(\mu,u,v): |\mu|\le 1, |u|\le 1, |v|\le 1\}$ corresponds to  a neighborhood of $P$ in $\Delta \times \C^2$ where $P$ has local coordinates $(0,0,0)$.  The coordinates may be chosen so that $\mu$ is the scaled parameter value and the sets with constant $u$ value, i.e. $\{(\mu, u_0, v)\}$, lie in a family of leaves of $\mathcal G^+_\lambda$.  In addition, we are able to choose the local coordinates so that the set with $v=0$, i.e. $\{(\mu, u, 0)\}$ corresponds to the family of local leaves of $J^-_\lambda$ through $P$.

The proof continues identical to that in \cite{BS6}, with $\ds \alpha: \mathcal J^-_+ \cap B \to \C$ by 
\[\alpha(\mu, u, v):=\frac{\varphi^+(\mu, u, v)}{\varphi^+(0,u,0)} \]
$\alpha$ is continuous with a discrete image, and it follows that $\varphi^+$ takes a common value on nearest neighbors. So $\varphi^+$ analytically continues to a neighborhood of $\mathcal J^-_+$. \BOX

\section{The solenoidal model}

\begin{center}
Why call it a ``solenoid?'' Pictured: $\pi_y J$ for $H(x,y)=(y,y^2-.3x)$.

\includegraphics[height=2in]{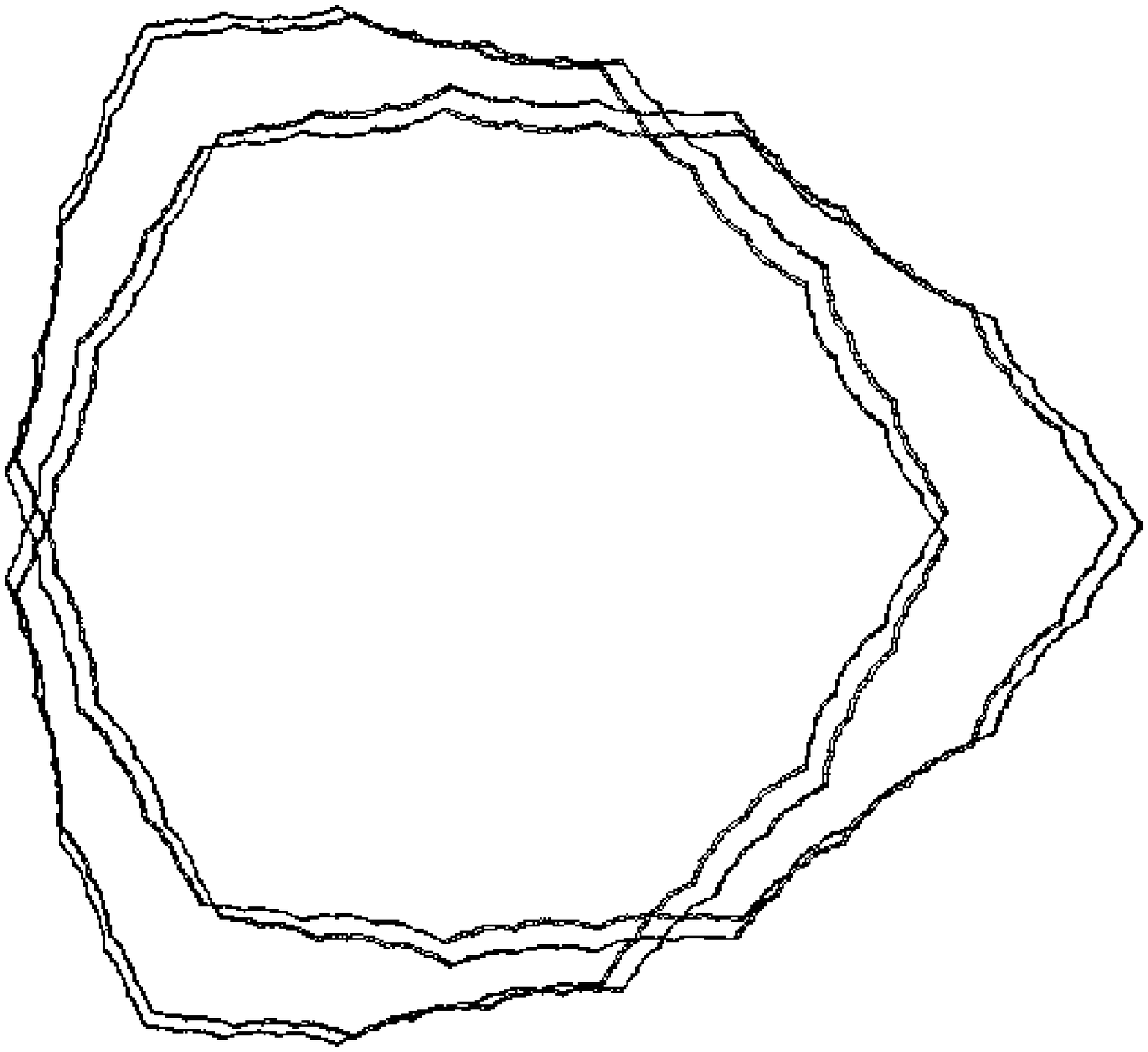}
\end{center}

\vspace{.1in}
\noindent Recall the following dynamical definitions (see \cite{S}):
\begin{defn} Given an invertible dynamical system $F:X \to X$, the {\bf unstable manifold} of $x\in X$,  $\ds W^u(x):= \{y\in X: \lim_{n \to \infty}  {\rm dist}(F^{-n}(x),F^{-n}(y))=0\}$. \\
\end{defn}
The following lemmas follow immediately from the definition.
\begin{lemma} $F(W^u(x))=W^u(F(x))$. \end{lemma}
\begin{lemma} $W^u(x)=W^u(y)$ if and only if $y\in W^u(x)$. \end{lemma}
\begin{defn} Let $F:X \to X$ be a dynamical system.  The {\bf inverse limit} of $X$ under $F$ is an invertible dynamical system, $\sigma: \hat X \to  \hat X$ where $\sigma$ is the shift map and  $\hat X$ is the space of bi-infinite sequences given by orbits of points in $X$ under $F$.  That is, $\hat X := \{ {\bf x} \in X^\mathbb Z: x_{t+1} = F(x_t) \; \;  \forall t \in \mathbb Z\} $ equipped with the product topology. \end{defn}

The degree-$d$ complex solenoid, $\Sigma$ is the inverse limit of $\C$ under $z \mapsto z^d$. The exterior solenoid, $\Sigma_+$ is given by restricting attention to $\C \setminus \mathbb D$.  Also let $\Sigma_0$ denote those ${\bf z} \in \Sigma$ with $|z_t|=1$. Let $s$ be a $(d-1)$th root of unity.  Consider the fixed points of $\Sigma_0$ given by $z_t = s$ for every $t\in \mathbb Z$.  Denote such a fixed point as ${\bf s}$.

\begin{prop} For the dynamical system given by $\sigma$ on $\ds \Sigma_0 \cup \Sigma_+$,\\ $\sigma \left (W^u(x) \right ) = W^u(x)$ if and only if $ W^u(x) = W^u({\bf s}) $ for some fixed point ${\bf s}$. \end{prop}

\pf{}
$\sigma \left (W^u({\bf s}) \right ) = W^u(\sigma({\bf s})) = W^u({\bf s})$, so one direction is clear.
Conversely, if $y \in W^u(x)$ and $y \in W^u(\sigma x) $, then ${\rm dist} (\sigma^{-n}(x), \sigma^{-n}(y))\to 0$ and ${\rm dist} (\sigma^{-n+1}(x), \sigma^{-n}(y))\to 0$.  This implies that ${\rm dist} (\sigma^{-n}(x), \sigma^{-n+1}(x)) \to 0$ as $n \to \infty$.  By the definition of distance given by the product topology on the solenoid, $|x_t-x_{t+1}| \to 0$ as $t\to -\infty$, i.e. $|x_t-x_t^d| \to 0$.  Hence, since $|x_t|\ge 1$, $x_t$ converges to a $(d-1)$th root of unity, $s$, as $t \to -\infty$.  Therefore $x \in  W^u({\bf s})$,  and it follows from Lemma 5.3 that $W^u(x) = W^u({\bf s})$. \BOX\\

Now we construct the semiconjugacy to the exterior solenoid.
Such a semiconjugacy requires a map $\Phi: J^-_+ \to \Sigma_+$ such that $\sigma \Phi = \Phi  H$.
Equivalently, require $\Phi_{t+1} =\Phi_t^d = (\Phi \circ H)_t$.
Now $\Phi_0:J^-_+ \to \C \setminus \mathbb D$ completely determines all coordinates $\Phi_t$ by the relation $\Phi_t=(\Phi \circ H^t)_0$, and $(\Phi \circ H(\lambda, p))_0= (\Phi_0(\lambda, p))^d$. Hence, $\varphi^+$ serves as this 0-th coordinate map.
Define the Bedford-Smillie solenoid map $\Phi: J^-_+ \to \Sigma_+$ by
\[\Phi_t := \varphi^+\circ H^t \;\;\;\;\; \text{for every }t\in \mathbb Z \]

\begin{thm} [\cite{BS6} Theorem 3.2, \cite{BS7} Theorem 4.3] If $H$ is hyperbolic and unstably connected, $\Phi$ is a continuous finite-degree covering map which is holomorphic and injective on the leaves of $J^-_+$ and such that \[\sigma \Phi = \Phi H\]
Under $\Phi$, the lamination of $J^-_+$ given by unstable manifolds is taken to the lamination (given by unstable manifolds) of $\Sigma_+$.\end{thm}

\begin{prop} If $p_0$ is a saddle fixed point of $H$ with $H$ hyperbolic and unstably connected, and $W^u(p_0)\cap U^+$ consists of only one component, then $\Phi(W^u(p_0) \cap U^+) = W^u({\bf s})\cap \Sigma_+$ for some $(d-1)$th root of unity $s$. \end{prop}
\pf{}
$\sigma \Phi(W^u(p_0) \cap U^+) = \Phi(H(W^u(p_0) \cap U^+)) = \Phi(W^u(p_0) \cap U^+)$. Proposition 5.5 and Theorem 5.6 imply $\Phi(W^u(p_0) \cap U^+) =  W^u({\bf s})\cap \Sigma_+$. \BOX

\begin{prop} If $s$ is a $(d-1)$-st root of unity and $q'\in J$ satisfies $\Phi(W^u(q') \cap U^+) = W^u({\bf s})\cap \Sigma_+$, then there exists $q \in W^u(q')$ such that $q$ is a periodic point of $H$, with period no more than $k$, the covering degree of $\Phi$.
\end{prop}

\pf{}  For $j\in \mathbb Z$, $W^u({\bf s})\cap \Sigma_+ =  \sigma^j(W^u({\bf s})\cap \Sigma_+)  = \sigma^j \Phi (W^u(q') \cap U^+) = \Phi (H^j(W^u(q') \cap U^+) )= \Phi (W^u(H^j(q')) \cap U^+)$.  Thus $W^u(H^j(q')) \cap U^+$ maps to $W^u({\bf s})\cap \Sigma_+$ under $\Phi$ for any $j$. $W^u({\bf s})\cap \Sigma_+$ has only $k$ preimages. So $W^u(q')$ is periodic, since $j$ may be negative as well as positive.  $W^u(q')=W^u(H^n(q'))=H^n(W^u(q'))$ for some $n \le k$.  There exists a conformal parametrization $\tau: \C \to W^u(q')$. Now $\tau^{-1}H^n\tau: \C \to \C$ is an expanding map (since $H$ is hyperbolic), so it has a fixed point $z$.  Set $q=\tau(z)$.  $H^n(q)=q$ and $W^u(q)=W^u(q')$. \BOX 

Let $\Theta: \Sigma_+ \to \Sigma_0$ by $\Theta({\bf z})_t = z_t/|z_t|$.  The ``external ray'' above $\omega \in \Sigma_0$ is the path having $\Theta \equiv \omega$ in $\Sigma_+$.  Hence, $J^-_+$ is equipped with solenoidal external rays via the fibers of rays in $\Sigma_+$ under the map $\Phi$. Two rays are said to be identified if they have the same landing point.

\section{The holomorphic motions respect the basic models}

As pointed out in \cite{M}, given a polynomial $P$ of degree $d$ with connected filled Julia set $K$, any bijective holomorphic map $\hat \phi:\C \setminus K \to \C \setminus \mathbb D$ such that $\hat \phi(z) \to \infty$ as $z\to \infty$ and satisfying  \[\hat \phi P(z)=  \big (\hat \phi(z)\big )^d \] is unique up to multiplication by a $(d-1)$-st root of unity.

The B\"ottcher coordinate is a holomorphic bijection $\phi: \C \setminus K \to \C \setminus \mathbb D$ satisfying $\phi P(z) = (\phi(z))^d$ and given by: \[\phi(z)=\lim_{n \to \infty} (P^n(z))^{1/d^n}\] on a neighborhood of infinity, where the roots are chosen to approximate the identity.  By the connectivity of $K$, the mapping $\phi$ uniquely extends to all of $\C \setminus K$.

We furnish an example to illustrate that the uniqueness of the B\"ottcher coordinate does not hold for quasiconformal maps.  First, interchange the role of $\infty$ with $0$ and consider $P(z)=z^2$.  Then $\hat \phi(z)=ze^{i \ln |z|}$ is a quasiconformal, non-holomorphic bijection defined on a neighborhood of $0$, with $\ds \lim_{z\to 0}\hat \phi(z)=0$ and satisfying the conjugacy $\hat \phi(z^2)=\big (\hat \phi(z)\big )^2$. 

If we assumed $\psi_\lambda$ to be holomorphic on $\C \setminus K_0$, then the proof of Theorem 1.1 would follow immediately from the uniqueness of the B\"ottcher coordinate as follows:  $\phi_\lambda \psi_\lambda$ satisfies the conjugacy, so by the uniqueness of the B\"ottcher coordinate, $\phi_\lambda \psi_\lambda=\zeta_\lambda \phi_0$, where  $\zeta_\lambda$ is a $(d-1)$th root of unity.  $\zeta_\lambda$ is  holomorphic in $\lambda$ and discrete-valued with $\zeta_0=1$. Hence $\phi_\lambda \psi_\lambda= \phi_0$.

However, since we only know $\psi_\lambda$ to be quasiconformal with harmonic Beltrami coefficient on $\C \setminus K_0$, we instead use the uniqueness of the Bers-Royden extension to prove the theorem.  It then follows that $\psi_\lambda$ must in fact be holomorphic on $\C \setminus K_0$ for each $\lambda\in \Delta$.

\pf{of Theorem 1.1} Given the McMullen-Sullivan holomorphic motion $\psi_\lambda$ defined on $\hat \C$ over a region of parameter space $\Delta$, for each $\lambda \in \Delta$ let
 \[\alpha_\lambda:= \phi_0^{-1}\phi_\lambda \psi_\lambda: \C \setminus K_0 \to \C \setminus K_0 \]

\begin{lemma}
$\alpha: \Delta \times \C \setminus K_0 \to \C \setminus K_0$ is a holomorphic motion such that $\alpha_\lambda$ has harmonic Beltrami coefficient on $\C \setminus K_0$ for each fixed $\lambda \in \Delta$.
\end{lemma}
\pf{of Lemma 6.1}
$\alpha_0(z)=\phi_0^{-1}\phi_0 \psi_0(z) = \psi_0(z)=z$. 
 $\alpha_\lambda$ is a homeomorphism for each fixed $\lambda$ since it is a composition of homeomorphisms and $\alpha_\lambda(z)$ is holomorphic in $\lambda$ for each fixed $z \in \C \setminus K_0$ since it is a composition of maps holomorphic in $\lambda$. As $\phi_0^{-1}\phi_\lambda$ is holomorphic for each fixed $\lambda\in \Delta$, the Beltrami coefficient remains harmonic. \BOX \\

By Theorem 3.2, $\alpha$ extends to a holomorphic motion on the closure, $\overline{ \C \setminus K_0}$ over $\Delta$.  Since  $\alpha_\lambda P_0 = P_0 \alpha_\lambda$, it follows that $\alpha_\lambda$ maps $J_0$ into itself for every $\lambda \in \Delta$.  Fixing any $z_0 \in J_0$, $\alpha_\lambda(z_0)$ is holomorphic as a function of $\lambda$. By the open mapping theorem of complex analysis, $\{ \alpha_\lambda(z_0) \}_{\lambda \in \Delta}$ is an open set or a constant. But since $\alpha_\lambda(z_0)$ must lie in $J_0$, which has empty interior, $\alpha_\lambda(z_0)$ must be a constant, namely $z_0=\alpha_0(z_0)$.  Since this holds for any $z_0 \in J_0$, $\alpha_\lambda$ is the identity map on $J_0$ for any $\lambda \in \Delta$.

By Theorem 3.4, the holomorphic motion $\alpha$ extends canonically to a motion of $\hat \C$ over $\Delta_{1/3}$.  This extension is also the canonical extension of $\alpha_\lambda|_{J_0}$, since the original $\alpha_\lambda$ has harmonic Beltrami coefficient in $\C \setminus K_0$. 
But the identity is also an extension of $\alpha_\lambda|_{J_0}$ with harmonic Beltrami coefficient.  By the uniqueness of the Bers-Royden theorem, $\alpha_\lambda$ must be the identity on all of $\overline{ \C \setminus K_0}$ for any $\lambda \in \Delta_{1/3}$.  Since $\alpha$ is holomorphic in $\lambda$, it follows that  $\alpha$ must be the identity on $\overline{ \C \setminus K_0}$ for all $\lambda \in \Delta$.

It immediately follows that \[ \phi_0(z)=\phi_\lambda \psi_\lambda(z) \] for all $\lambda \in \Delta$ and $z\in \C \setminus K_0$. Furthermore, the McMullen-Sullivan holomorphic motion $\psi_\lambda=\phi_\lambda^{-1}\phi_0$ is holomorphic on $\C \setminus K_0$ for every $\lambda \in \Delta$. This completes the proof of Theorem 1.1. \BOX
\vspace{.2in}

Turning to higher dimensions, again, since we only know $\Psi_\lambda$ to be leaf-wise quasiconformal, not necessarily holomorphic, we cannot use the uniqueness of $\varphi^+$ (as explained in \cite{HO1}, page 40) to give a short proof.  Instead, using the theorem of Bers and Royden, it follows that the solenoidal mapping $\Phi$ commutes with the holomorphic motion $\Psi$.  A surprising consequence is the fact that the Buzzard-Verma motion $\Psi$ is holomorphic on each leaf of $J^-_+(0)$.

\pf{of Theorem 1.2} Let $H_0$ be hyperbolic and unstably connected. From \cite{O}, there exists a saddle fixed point $p_0$ of $H_0$ with $W^u(p_0)\cap U_0^+$ consisting of a single, simply connected component. Since \[H_\lambda \Psi_\lambda(p_0) = \Psi_\lambda(H_0(p_0)) = \Psi_\lambda(p_0)\]  $\Psi_\lambda(p_0)$ is a fixed point of $H_\lambda$.   $\Psi_\lambda(p_0)$ is a saddle fixed point since $\Psi_\lambda$ preserves the unstable and stable manifold splitting. $\Psi_\lambda(W^u(p_0)\cap U_0^+) = W^u(p_\lambda)\cap U_\lambda^+$ must also consist of a single, simply connected component ($\Psi_\lambda$ is a homeomorphism).

By Proposition 5.7, let $\Phi_0^{-1}$ denote the branch of the inverse of the covering map $\Phi_0$ such that $\Phi_0^{-1}(W^u({\bf s}) \cap \Sigma_+) = W^u(p_0)\cap U^+_0$, where ${\bf s}$ is a fixed point of the solenoid.
Let \[f_\lambda:= \Phi_0^{-1}\Phi_\lambda\Psi_\lambda: W^u(p_0)\cap U^+_0 \to W^u(p_0)\cap U^+_0\]

Let $E:= \tau^{-1}_0( W^u(p_0)\cap U^+_0)$ and let  \[F_\lambda:= \tau^{-1}_0 f_\lambda \tau_0:E \to E \]
Like Lemma 6.1, $F_\lambda$ is a holomorphic motion of $E$. For any $z\in E$,  
$F_0(z)=\tau^{-1}_0 \Phi_0^{-1}\Phi_0 \Psi_0 \tau_0 (z) = z$ by the definition of $\Phi_0^{-1}$.
 $F_\lambda$ is a homeomorphism for each fixed $\lambda$ since it is a composition of homeomorphisms and $F_\lambda(z)$ is holomorphic in $\lambda$ for each fixed $z \in E$ since it is a composition of maps holomorphic in $\lambda$. Furthermore, $F_\lambda$ has harmonic Beltrami coefficient on $E$ for each fixed $\lambda \in \Delta$.

Again by Theorem 3.2, the holomorphic motion $F_\lambda$ extends to the closure of $E$ in $\C$. 
Pulling back to the unstable manifold, the motion $f_\lambda$ extends to the closure of $W^u(p_0)\cap U^+_0$ in $W^u(p_0)$, that is, to $W^u(p_0)\cap J^+_0 =W^u(p_0)\cap J_0$ via $f_\lambda:=\tau F_\lambda \tau^{-1}$.

Restricting to this boundary, it follows that $\tau^{-1} (W^u(p_0)\cap J_0)$ has empty interior.
Following the proof of Theorem 1.1, $F_\lambda$ is the identity on the boundary of $E$ and $f_\lambda$ must be the identity on $W^u(p_0)\cap J_0$ .  Furthermore, by the uniqueness of the Bers-Royden extension, $F_\lambda$ is the identity on $\C$ and hence  $f_\lambda$ is the identity map on $W^u(p_0)$.
\[ \Phi_\lambda \Psi_\lambda=\Phi_0: W^u(p_o)\cap U^+_0 \to W^u({\bf s}) \cap \Sigma_+ \]
Now $W^u(p_0) \cap U^+_0$ is dense in $J^-_+(0)$.  For each $\lambda\in \Delta$, $\Phi_\lambda \Psi_\lambda$ is continuous on $J^-_+(0)$ and $\Phi_0$ is also continuous on $J^-_+(0)$.   These continuous maps agree on a dense subset, therefore 
\[ \Phi_\lambda \Psi_\lambda=\Phi_0:J^-_+(0) \to \Sigma_+ \]
 Furthermore, the Buzzard-Verma holomorphic motion $\Psi_\lambda$ is holomorphic on each leaf of $J^-_+(0)$. This completes the proof of Theorem 1.2. \BOX

Corollary 1.3 immediately follows.  In general, we see that $\Phi$ varies nicely within regions of hyperbolic unstably connectedness.  In the next section we prove Theorem 1.4, in which case $\Phi_\lambda$ has a constant covering degree of one, that is, $\Phi_\lambda$ is a homeomorphism.

\section{For small Jacobian}
As was mentioned before its statement, the hypothesis of Theorem 1.4 is not vacuous.
From \cite{HO2}, a slight perturbation, $H_b$ for $b$ near $0$,  of a hyperbolic polynomial $P$ is also hyperbolic. Assuming the Julia set of $P$ is connected, the inverse limit of the Julia set of $P$ is connected.  Hence by the homeomorphism given in \cite{HO2} from the inverse limit of the Julia set to that of the perturbed H\'enon map $H_b$, the Julia set of $H_b$ is connected as well.  By \cite{BS6} Theorem 0.2, this means $H_b$ is unstably connected.

We prove Theorem 1.4 by restricting attention to external rays that land at periodic points. The linearizing coordinate near a periodic point is locally injective, giving a neighborhood in the unstable manifold on which we locally extend the Buzzard-Verma motion in the case when $b=0$.  The external ray continues this extension to a larger neighborhood that intersects $\mathcal V^+$.  As before, the extended Buzzard-Verma motion still follows level sets of $\varphi^+$.  This gives an extension of $\varphi^+_0$ which must agree with the B\"ottcher coordinate.  Hence the external rays are preserved at $b=0$. Since the B\"ottcher coordinate is injective, we conclude that the solenoid map is injective.

When the Jacobian $b=0$, $H_0$ collapses $\C^2$ onto the parabola $\mathcal W:=\{(x,y)\in \C^2: y=P(x)\}$ and $H_0$ is not a diffeomorphism.  Since $H_0(x,y)=(y, P(y))$ the dynamics are completely described by the one-dimensional polynomial map $P$.

Let $\Delta$ denote the set $\{|b|<\epsilon\} \subset \C$.  Write $\Delta ^* := \Delta \setminus \{0\}$. Assume $\Delta ^*$ lies in $\mathcal M_P$. Throughout this section we will shrink $\epsilon>0$ as needed. By showing the identification of rays and the injectivity of $\Phi$ on $\Delta^*$ for some $\epsilon$, Corollary 1.3 implies the identification of rays and the injectivity of $\Phi$ on $\mathcal M_P$.

Suppose $\theta_0, \hdots , \theta_{K-1}$ are the angles of a periodic external ray for $P$ where  $\theta_{j+1}=d \theta_j \; ({\rm mod} \; 2\pi)$ for $j=0, \hdots, K$ with $\theta_K=\theta_0$. A corresponding solenoidal external ray is one that maps under $\Phi_b$ to the same periodic cycle of angles.  In fact, it will be shown that such a ray is unique and

\begin{thm}\label{rays} If $b \in \mathcal M_P$, then the identifications of periodic external rays are the same for $P$ and $H_b$.   \end{thm}

\pf{} 
First we extend the Buzzard-Verma holomorphic motion to the case when $b=0$ in a neighborhood of an external ray over a periodic point $q_0$.

 Let $q_0' \in \C$ be a periodic point in the Julia set of $P$ with period $m$.  Write $q_k':=P^k(q_0')$ for $k=0, \hdots, m$, with $q_0'=q_m'$. $\ds M:= d^m\Pi _{k=1}^m q_k'$ is the multiplier of the periodic orbit of $q_0'$. Since $P$ is hyperbolic, $|M|>1$.  Write $q_0:=(q_0', q_{1}')\in \C^2$. So $H_0^m(q_0)=q_0$.

Let $R_0$ be an external ray in $\C$ that lands at $q_0'$. Let $R_j=P^j(R_0)$.
Each $R_j$ is an external ray in $\C$ that lands at $q_k'$, $k=(j\; {\rm mod}\; m)$. There may be multiple, ``identified'' rays that land at the same (pinch) point, and each ray must be periodic of period $mN$ for some $N<\infty$.  By definition, $R_j=\phi^{-1}(\{re^{i\theta_j}\})$ where $r>1$ and $\theta_{j-1}d\equiv \theta_j \;( {\rm mod}\; 2\pi)$ for $j=1, \hdots, mN$. 

Let $\mathcal W_U = \mathcal W \cap U^+_0$. Let $\Gamma_j:=\{(x,y)\in \mathcal W_U: x\in R_{j-1}, y\in R_j\}$. Let \[T_{j}:= \big \{(x,y)\in \mathcal W_U: |\arg \phi(x) - \theta_{j-1}|<\frac{\pi}{d+1}, |\arg \phi(y) - \theta_j|<\frac{\pi}{d+1}\big \}\] \[T := \bigcup_{j=1}^{mN} T_{j}\]

Suppose $(x,y) \in T_j$ for some $j$. Suppose $\phi(x)=re^{i(\theta_{j-1}+\zeta)}$ and $\phi(y)=se^{i(\theta_{j}+\xi)}$, where $|\zeta|<\frac{\pi}{d+1}$, $|\xi|<\frac{\pi}{d+1}$, $r>1$, and $s>1$. Since $y=P(x)$, $r=s^{1/d}$ and $(\theta_{j-1}+\zeta)d\equiv \theta_j+\xi \;( {\rm mod}\; 2\pi)$. The only solution is $\zeta =\xi/d$.  As such, we may define a branch of $H_0^{-1}: T \to T$ by $H_0^{-1}(y,P(y)) = (x,y)$ where $x$ is completely determined by $y$. Furthermore, this also implies $\pi_y$ is injective on $T_j$, where $\pi_y: \C^2 \to \C$ denotes projection onto the second coordinate. 

  Let $F: \Delta \times \C^2 \to \C^2$ be given by $F(b,q)= H^m_b(q)-q$.  Observe,  \[DH_0^m(q_0) = \left [\begin{array}{cc} 0 & M/(dq_1')\\ 0 & M \end{array}\right ] \] has eigenvalues $0$ and $M$.
Hence $\ds \det \frac{\partial F}{\partial q}(0, q_0)=M-1 \ne 0$.  So the implicit function theorem applies to $F$: For some $\epsilon>0$ there is a holomorphic function $q(b)$, with $q(0)=q_0$ and $q(b)$ is the unique corresponding period-$m$ point for $H_b$, with $|b|<\epsilon$. Write $q_b$ for $q(b)$.

Notice for simply connected neighborhoods in $\Delta^*$, by the uniqueness of the implicit function theorem, $q(b)$ must agree with the Buzzard-Verma motion defined here.  That is, for $b_0,b_1 \in \Delta^*$, and the motion $\Psi$ centered at $b_0$, \[ \Psi_{b_1}(q_{b_0})=q_{b_1}\]

Turning attention to the unstable manifolds of $q_b$ and following \cite{HP}, consider the linearizing coordinate $\gamma_b:\C \to \C$ for each $b \in \Delta$ (with $\gamma_b(\C)=W^u(q_b)$ for $b\ne 0$) given by \[\gamma_{b,n}(z):= H_b^{mn}\left (q_b+\frac{z}{\lambda_b^n}v_b \right )\]
\[\gamma_b(z):=\lim_{n \to \infty} \gamma_{b, n}(z)\]
where $v_b$ is an eigenvector in $\C^2$ of $DH^m_b(q_b)$ with eigenvalue $\lambda_b$ such that $|\lambda_b|>1$.  Observe when $b=0$, $\lambda_0=M$ and $v_0=\left [ \begin{array}{c}1/2q_1'\\ 1 \end{array} \right ] $. The parametrization $\gamma$ is such that \[H^m_b(\gamma_b(z))=\gamma_b(\lambda_b z)\] (hence ``linearizing coordinate'') and $\gamma_b(0)=q_b$ for every $b \in \Delta$.

We may choose eigenvectors so that $\gamma_{b,n}(z)$ is holomorphic in $b$ and in $z$.  From Hubbard and Papadantonakis' proof of convergence in the definition of $\gamma$ (see \cite{HP}), it is apparent that $\gamma$ is the locally uniform limit of holomorphic functions; therefore, the linearizing coordinate $\gamma$ is holomorphic in $b$ and $z$ on $\Delta \times \C$.

  Now $\gamma_0$ in general is not injective.  (For example, in the simple case where $P(y)=y^2$, at the fixed point $q_0=(1,1)$, we obtain $\gamma_0(z)=(e^{z/2}, e^z)$, an infinite-to-one mapping). 
 However, differentiating the functional equation above yields \[ \frac{\partial \gamma_b}{\partial z}(0)= v_b\]  $\gamma_0$ must be locally injective on a small neighborhood $N$ of $0$.  By shrinking $N$ and $\epsilon$, it follows that $\gamma_b$ is injective on $N$ for any $b\in \Delta$.

Recall that $G_b^+(p)$ is pluriharmonic on $\mathcal U^+$, meaning it is locally the real part of a holomorphic function on $\mathcal U^+$.  The level sets of this holomorphic function correspond to the plaques of $\mathcal G_b^+$, hence $\mathcal G_b^+$ is a holomorphic family of laminations (see \cite{BV}).

 $\gamma_b(N)$ is transverse to $\mathcal G_b^+$ for $b\in \Delta^*$ since the unstable manifolds are transverse to  $\mathcal G_b^+$ in the unstably connected case, as in \cite{BV}. $\gamma_0(N)$ is also transverse to $\mathcal G_0^+$ by observation: The non-transverse points must be given by critical points of $P$, which lie in the interior of $K^+$, away from the periodic point $q_0$. Hence we may assume that $\gamma_0(N)$ lies on the parabola $\mathcal W$ away from any points of non-transversality with $\mathcal G_b^+$.

Let \[\hat N_U=U^+_0 \cap \gamma_0(N)\;\;\; {\rm and }\;\;\; N_U=\gamma_0^{-1}(\hat N_U) \]


Note that in general $\hat N_U$ may consist of multiple components, corresponding to the number of external rays that land at $q_0$. Now $\gamma: \Delta \times N \to \C^2$ is a holomorphic map such that for each fixed $b\in \Delta$, $\gamma_b$ is injective and $\gamma_b(N)$ is transverse to $\mathcal G^+_b$ in $U^+_b$.  By Lemma 4.1 of Buzzard and Verma \cite{BV}, as a consequence of the definition of lamination families and the implicit function theorem, for some $\epsilon>0$ there exists a function $\Upsilon: \Delta \times \hat N_U  \to \C^2$ with the following properties: For each $p_0 \in \hat N_U$, $\Upsilon_0 (p_0)=p_0$, $\Upsilon_b(p_0)$ is holomorphic in $b\in \Delta$, and $\Upsilon_b(p_0)$ is given by the intersection of the family of leaves of $\mathcal G_b^+$ through $p_0$ with the subset of the unstable manifold of $q_b$ that lies in $\gamma_b(N)$. 

This map $\Upsilon$ must be a local extension of the Buzzard-Verma motion $\Psi$. Let $b_0,b_1 \in \Delta^*$ and let $\Psi$ be the Buzzard-Verma holomorphic motion centered at $b_0$, i.e. $\Psi_{b_0}$ is the identity. The following diagram commutes:
\[ \xymatrix{
 & \hat N_U \ar[dl]_{\Upsilon_{b_0}} \ar[dr]^{\Upsilon_{b_1}} & \\
J^-_+(b_0) \ar[rr]_{\Psi_{b_1}} & & J^-_+(b_1) } \]
that is, for $p \in \hat N_U$, we have $\Psi_{b_1}\Upsilon_{b_0}(p)= \Upsilon_{b_1}(p)$ since the definition of $\Psi$ for $b_1\in \Delta^*$ in \cite{BV} Theorem 5.7 agrees with the definition of $\Upsilon$ here. More precisely, $\Psi_{b_1}\Upsilon_{b_0}(p)$ is given by the intersection of the family of leaves of $\mathcal G_b^+$ through $\Upsilon_{b_0}(p)$ with the unstable manifold of $q_{b_1}$.   $\Upsilon_{b_1}(p)$ is given by the intersection of the family of leaves of $\mathcal G_b^+$ through $p$ with the unstable manifold of $q_{b_1}$.  Since the family of leaves through $p$ is the same as the family of leaves through $\Upsilon_{b_0}(p)$, these two motions agree.  Furthermore, since $\Upsilon_0$ is injective, it follows that $\Upsilon_{b_0}$ must also be injective.

For fixed $b_0 \in \Delta^*$, we have $\Psi_{b_1}=\Upsilon_{b_1} \Upsilon^{-1}_{b_0}$ on $\Upsilon_{b_0}(\hat N_U)$.  For $b_1=0$, $\Upsilon_{b_1}$ is the identity. Letting $b_1 \to 0$, we define \[ \Psi_0|_{\Upsilon_{b_0}(\hat N_U)} =\Upsilon_{b_0}^{-1}\]  So our motion $\Psi$ locally extends to $b=0$. Theorem 1.2 shows $\Psi_{b_1}$ is holomorphic on leaves of $J^-_+(b_0)$.  So the limit ($b_1 \to 0$) function, $\Psi_0|_{\Upsilon_{b_0}(\hat N_U)}$ is holomorphic and its inverse, $\Upsilon_{b_0}$, is holomorphic on $\hat N_U$.
For $p\in U^+(b_0)$:
\[\Psi_{b_1}H_{b_0}^m(p)=H_{b_1}^m\Psi_{b_1}(p)\] 
so for $p\in \Upsilon_{b_0}(\hat N_U)$ such that $H^m_{b_0}(p) \in \Upsilon_{b_0}(\hat N_U)$:
\[\Upsilon_{b_1} \Upsilon_{b_0}^{-1} H_{b_0}^m(p)= H_{b_1}^m \Upsilon_{b_1}\Upsilon_{b_0}^{-1}(p)\]
Now let $b_0 \to 0$. By continuity, 
\[\Upsilon_{b_1}  H_0^m(p)= H_{b_1}^m\Upsilon_{b_1}(p)\]

\includegraphics[height=2.75in]{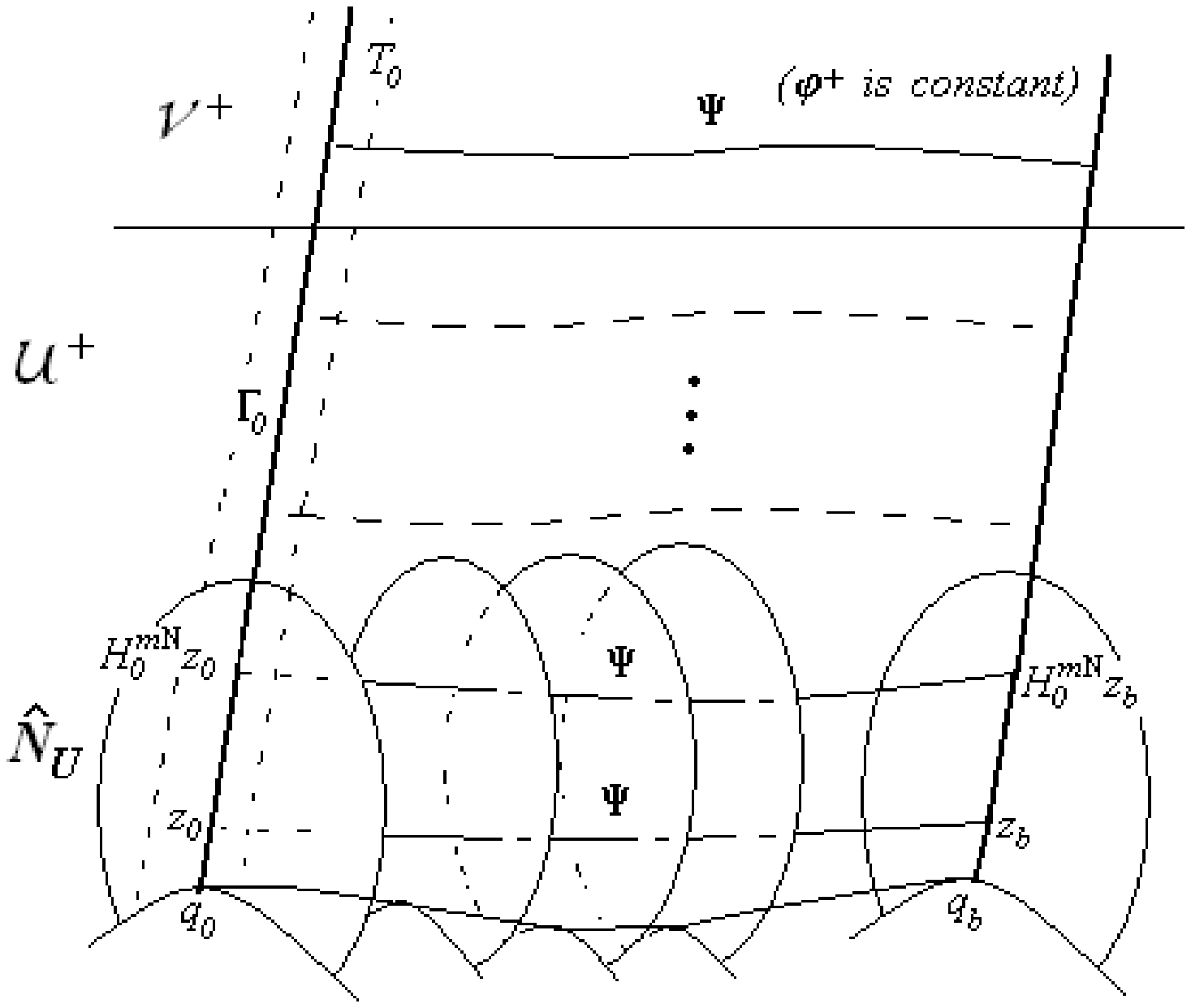}

Next, we want to extend $\Upsilon_b$ to a larger domain.
Define $\Upsilon_b$ on $T$ as follows: \[ \Upsilon_b(z) = H_b^{jmN}(\Upsilon_b(z_0))\] where $z_0:=H_b^{-jmN}(z) \in \hat N_U$ for some $j$. 
This is well-defined by the remarks above, still a holomorphic motion, and maintains the conjugacy $H_b^m\Upsilon_b=\Upsilon_bH_0^m$.

With $\Upsilon_b$ defined on $T$, we now use Theorem 1.2 to extend the definition of $\varphi^+_0$ to $T$.
For $b\in \Delta^*$, $\varphi^+$ is holomorphic on a neighborhood of $\mathcal J^-_+$ by Lemma 4.6,  and for fixed $p\in J^-_+(0)$, $\varphi^+_b \Upsilon_b(p)$ is constant as a function of $b \in \Delta^*$. 
Define $\ds \hat \varphi^+_0: T \cup \hat N_U \to \C \setminus \mathbb D$ by \[\hat \varphi^+_0(p):= \varphi^+_{b_0} \Upsilon_{b_0}(p)\]
for any $b_0\in \Delta^*$; the choice of $b_0$ is irrelevant, as pointed out above.  Since $\Upsilon_{b_0}$ is holomorphic, $\hat \varphi^+_0$ is also holomorphic. Recall that $\varphi^+_0$ is already defined on $V_0^+$. By continuity, $\varphi_0^+ = \hat \varphi^+_0$ on $V_0^+ \cap T$. Rename $\hat \varphi^+_0$ as simply  $\varphi_0^+$.

From the definition of the B\"ottcher coordinate, $\varphi^+_0(p)=\phi(\pi_y p)$ for $p\in V^+_0$. Earlier we showed $\pi_y$ is injective on each $T_j$. So $\varphi^+_0 \pi_y^{-1}$ and $\phi$ are holomorphic on $\pi_y(T_j)\subset \C$ and agree on the subset $\pi_y(V_0^+ \cap T_j)$.  Hence  $\varphi^+_0 \pi_y^{-1} = \phi$ on $\pi_y(T_j)$ and  $\varphi^+_0 = \phi \pi_y$ on $T_j$.

Hence, for any periodic external ray $R_j$ of $P$ that lands at $q_k'$ with corresponding ray $\Gamma_j \in \mathcal W$, \[\varphi^+_b \Upsilon_b(\Gamma_j)=\varphi^+_0(\Gamma_j)=\phi \pi_y(\Gamma_j)=\phi(R_j)\] has constant argument, i.e. $\Upsilon_b(\Gamma_j)$ is a periodic external ray for the map $H_b$ which lands at $q_b$.

Next we show the converse direction: that any periodic ray landing at $q_b$ has a unique corresponding ray landing at $q_0'$. This will follow from a counting argument.  
Notice that $\gamma_b^{-1}\Upsilon_b \gamma_0$ forms a holomorphic motion of $N_U$. By Theorem \ref{mss} this motion extends to the closure and hence defines $\Upsilon_b$ on $\gamma_0(N)\cap J_0$.  By continuity, the extended $\Upsilon$ is still simply an extension of $\Psi$ to $b=0$.
Essentially, $\Upsilon_b$ is a homeomorphism from a neighborhood of $q_0$ in $\gamma_0(N)\cap \overline{U_0^+}$ to a neighborhood of $q_b$ in $\gamma_b(N)\cap \overline{U_b^+}$.  What is not immediate is that the components of $\Upsilon_b (\hat N_U)$ coincide with those for $\gamma_b(N)\cap U_b^+$.  We will bound how far points in $J_0$ can move to ensure that it is impossible for some overlooked component to creep into $\gamma_b(N)\cap U_b^+$.

Let $B_{\eta}(p)$ denote the ball centered at $p$ with radius $\eta$ in $\C^2$.  Choose $\eta>0$ small such that $\mathcal W \cap B_{4\eta}(q_0) \subset \gamma_0(N)$.  Let $M\subset N$ be a neighborhood of 0.  We may suppose that $\epsilon$ and $M$ are small enough that $\gamma_b(M)\subset B_{\eta}(q_0)$. 

To bound how far our points can move, we consider a uniform cone field on a neighborhood of $\mathcal W \cap V(R_\Delta)$.
 For each point $p \in \mathcal W\cap V(R_\Delta)$ and $v\in \C^2$, let $\|v\|^s_p$ be the length of $v$ in the stable direction (i.e. $x$) and let $\|v\|^u_p$ be the length of $v$ in the unstable direction (i.e. tangent to $\mathcal W$). 
These definitions can be extended to include $p\in\mathcal W'$, a small neighborhood of $\mathcal W \cap V(R_\Delta)$ in $\C^2$. Let $C_p^s=\{v\in \C^2: \|v\|^s_p \ge \|v\|^u_p\}$. Choose $\mathcal W'$ small enough that for any $p\in \mathcal W'$ and $v\in C^s_p$, if $p+v \in \mathcal W'$ then $\|v\|<\eta$ in $\C^2$.
We may choose $\epsilon$ small enough that $J_b\subset \mathcal W'$ and if $p_b\in J_b$, then its local stable manifold $W^s_{loc}(p_b)$ lies inside the cone $C^s_{p_b}$.

Let $\pi_x^{-1}: \C \to \mathcal W$ be given by $\pi_x^{-1}(z)=(z, P(z))$.  
Let Per$_0$ denote the set of periodic points of $H_0$. Since $P$ is hyperbolic, any periodic point of $H_0$ can be followed by the implicit function theorem at $b=0$ and continued by $\Psi$, as pointed out above. 
Given $p_0 \in$Per$_0$, let $g_b(p_0)$ denote the intersection of $\mathcal W$ with $W^s_{loc}(p_b)$, the local stable manifold of the corresponding point $p_b$.
Now define $a_b$ on $\pi_x($Per$_0)\subset \C$ by $a_b:=\pi_x g_b \pi_x^{-1}$.  
$a$ is a holomorphic motion as it is derived from a leafwise trivial holomorphic family of laminations.
By Corollary 2 of \cite{BR}, we can suppose $\epsilon$ is small enough that $\|g_b(p_0)-p_0\|<\eta$ for any $|b|<\epsilon$.  
Since $g_b(p_0)-p_b$ lies in the cone $C^s_{p_b}$, we have $\|g_b(p_0)-p_b\|<\eta$.  Also, if $w\in \mathcal W' \cap W^s_{loc}(p_b)$ then $w-p_b$ also lies in the cone $C^s_{p_b}$, and we have $\|w-p_b\|<\eta$. Hence, $\|w-p_0\|<3\eta$.

Let $\mathcal O$ be a component of $\gamma_b(M)\cap U_b^+$ such that $q_b \in \partial \mathcal O$. By local product structure and Proposition 5.7 of \cite{BS7}, $\partial \mathcal O$ must contain a point $w\ne q_b$ that lies in the local stable manifold of some periodic point $p_b\in J_b$.
As $w\in \gamma_b(M)$, $\|w-q_0\|<\eta$ which implies $\|q_0-p_0\|<4\eta$, i.e. $p_0\in \gamma_0(N)\cap J_0$.
Hence, $\Upsilon_b(p_0)=w$ and in particular, $\partial \mathcal O \setminus \{q_b\}$ contains points in $\Upsilon_b(\gamma_0(N)\cap J_0)$.

Consequently, as $\Upsilon_b$ is a homeomorphism, the number of components of $\gamma_b(M) \cap U_b^+$ that contain $q_b$ in their boundary must equal the number of components of $\gamma_0(M) \cap U_0^+$ that contain $q_0$ in their boundary.  The number of external rays in $W^u(q_b)\cap U_b^+$ that land at $q_b$ is equal to the number of components of $W^u(q_b)\cap U_b^+$ (A consequence of \cite{BS7} detailed in section 2.4.1 of \cite{O}). From one-dimensional dynamics, the number of components of $\gamma_0(M) \cap U_0^+$ with $q_0$ in their boundary is equal to the number of rays in $\C$ that land at $q_0'$.  Hence, there cannot be more rays landing at $q_b$ than the number of rays that land at $q_0'$.

Therefore, the periodic external rays in $W^u(q_b) \cap U_b^+$ that land at $q_b$ are in one-to-one correspondence with the external rays in $\C$ that land at $q_0'$. \BOX

\pf{of Theorem 1.4} By Theorem 5.6 and Proposition 5.8, any external ray $\Gamma$ of angle $0$ (i.e. a fiber under $\Phi_{b}$ of the ray over the fixed point ${\bf 1}$ in $W^u({\bf 1})\cap \Sigma_+$) must be a periodic external ray that lands (by \cite{BS7} 2.1) at a periodic point $q_{b}$ for $H_{b}$ (of period no more than $k$).   However, by Theorem \ref{rays}, the periodic external rays given by $\Phi_{b}$ for $H_{b}$ correspond to the periodic external rays given by $\phi$ for $P$.  Since $\phi$ is injective, it follows that there is only one such solenoidal ray landing at ${\bf 1}$. Consequently, $\Phi_{b}$ is injective. \BOX

Notice that in these cases where $\Phi$ is a homeomorphism, $\Phi^{-1}$ coincides with the homeomorphism known as $\Psi$ (not to be confused with the Buzzard-Verma motion) in Bedford and Smillie's work on solenoidal external rays (Theorem 4.1 of \cite{BS7}).

Empirical data seems to indicate the existence of hyperbolic, unstably connected H\'enon maps that are not simply a perturbation of one-dimensional hyperbolic polynomials.  For example, \cite{O} conjectures this to be the case for $H(x,y)=(y, y^2 - 1.124 - 0.125x)$.  The injectivity of $\Phi$ remains to be seen in these cases.

\begin{center}
\includegraphics[height=2.5in]{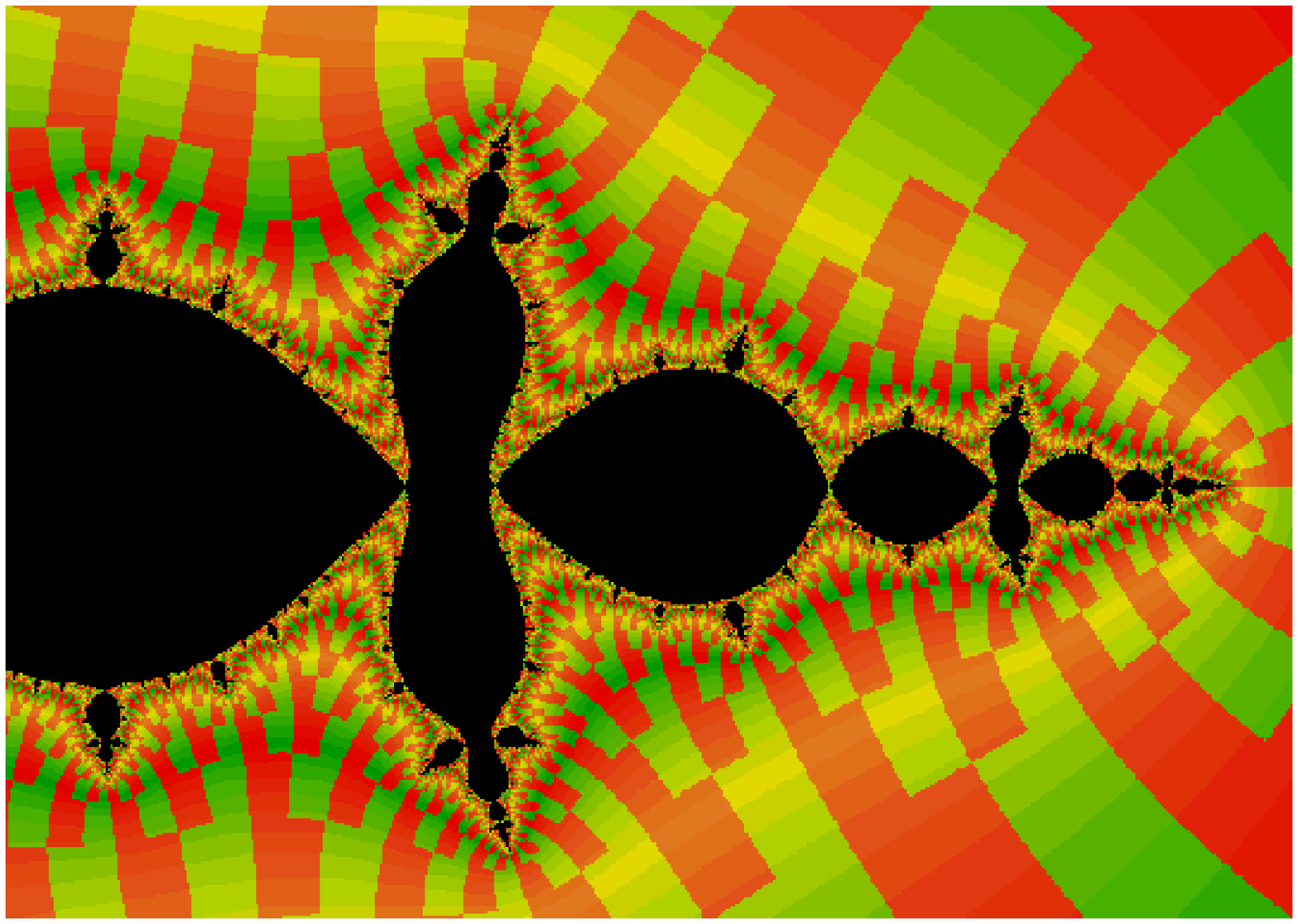}

The unstable manifold of the saddle fixed point for $H(x,y)=(y, y^2-1.124-.125x)$

 Picture created with FractalAsm: http://www.math.cornell.edu/$\sim$dynamics/
\end{center}

\pagebreak

\end{document}